\documentclass[a4paper, 12pt]{amsart}

\usepackage{amssymb,amsthm,amstext,amsmath,amscd,latexsym,amsfonts,amscd}
\usepackage{color}

\def\M{\mathcal M}
\def\S{\mathcal S}

\def\p{\mathfrak p}
\def\n{\mathfrak n}
\def\q{\mathfrak q}
\def\m{\mathfrak m}

\def\Ass{\mathrm {Ass}}
\def\Min{\mathrm {Min}}
\def\depth{\mathrm{depth}}
\def\Hom{\mathrm {Hom}}
\def\injdim{\mathrm{inj\,dim}}
\def\Spec{\mathrm {Spec}}
\def\Supp{\mathrm {Supp}}

\def\RMod{R\text{-}\mathrm{Mod}}
\def\Rmod{R\text{-}\mathrm{mod}}

\def\G{\varGamma }
\def\FS{\M_{f.s.} }
\def\FG{\S_{f.g.} }
\def\AR{\S_{Artin } }
\theoremstyle{plain} 
\newtheorem{thm}{\textbf Theorem}[section]

\newtheorem{prop}[thm]{\textbf Proposition}

\theoremstyle{definition}
\newtheorem{df}[thm]{\textbf Definition}
\newtheorem{rem}[thm]{\textbf Remark}
\newtheorem{exm}[thm]{\textbf Example}

\title[An example of Melkersson subcategory]{An example of Melkersson subcategory which is not closed under injective hulls}
\author{Takeshi Yoshizawa}
\keywords{Melkersson subcategory, Section functor, Serre subcategory}
\subjclass[2010]{13C60, 13D45}
\begin{document}
\begin{abstract}
The Melkersson subcategory is a special Serre subcategory which satisfies useful conditions $C_{I}$ defined in \cite{AM-2008}. 
It was proved that a Serre subcategory which is closed under injective hulls is a Melkersson subcategory.
However, it has been an open question whether the contrary implication holds. 
In this paper, we shall show that this question has a negative answer in general. 
\end{abstract}
\maketitle
\thispagestyle{empty}
\section{Introduction}
Throughout this paper, all rings are commutative noetherian ring, all modules are unitary and $R$ denotes a ring.  
We assume that all full subcategories $\S$ of the modules category $\RMod$ and the finitely generated $R$-modules category $\Rmod$ are closed under isomorphisms, that is if $M$ is in $\S$ and $R$-module $N$ is isomorphic to $M$ then $N$ is in $\S$.

\vspace{3pt}
In \cite{AM-2008}, 
M.\ Aghapournahr and L.\ Melkersson gave a useful condition $C_{I}$ on the Serre subcategory $\S$ of $\RMod$ where $I$ is an ideal of $R$.    
It is said that $\S$ satisfies the condition $C_{I}$ if the following condition holds: 
if $M=\G_{I}(M)$ and $(0:_{M}I)$ is in $\S$, then $M$ is in $\S$. 
They showed that local cohomology modules and Serre subcategories which satisfy such a condition have affinity for each other. 
After of this, 
the Serre subcategory which satisfies the condition $C_{I}$ for all ideals $I$ of $R$ was named Melkersson subcategory 
by M.\ Aghapournahr, A.\ J.\ Taherizadeh and A.\ Vahidi in \cite{ATV-2009}. 
For example, all Serre subcategories which are closed under injective hulls are Melkersson subcategory. 
So it is natural to ask the following question which was given in \cite{AM-2008}:   

\vspace{7pt}
\noindent 
{\bf Question.}~  
Is Melkersson subcategory closed under injective hulls?

\vspace{12pt}
In this paper, 
we shall show that this question has a negative answer in general. 
To be more precise, 
we denote by $\FG$ the Serre subcategory of all finitely generated $R$-modules 
and by $\FS$ the Serre subcategory of all $R$-modules with finite support. 
We shall see that a class 
\[ 
(\FG, \FS) = \left\{ X \in \RMod \text{ \Large $\mid$ }  
\begin{matrix} \text{there are $S \in \FG$ and $M \in \FS$ such that} \cr 
\minCDarrowwidth1pc \begin{CD}0 @>>> S @>>> X @>>> M @>>> 0\end{CD} \ \ \text{is exact.}\cr \end{matrix} \right\} 
\]
is Melkersson subcategory which is not closed under injective hulls on the ring of formal power series $R=k[[x, y]]$ in the indeterminate $x$ and $y$ with the coefficients in a field $k$.

\vspace{5pt}
The organization of this paper is as follows.
 
In section 2, 
we shall recall definitions of Melkersson subcategory and classes $(\S_{1}, \S_{2})$ of extension modules of a Serre subcategory $\S_{1}$ by another Serre subcategory $\S_{2}$. 
In section 3, we shall give a proof of main result. 
In Section 4, we shall see several remarks on Melkersson subcategory.

\vspace{7pt}
\section{Preliminaries}
In this section, 
we shall recall  several definitions which are necessary to prove the main result of this paper.

\vspace{5pt}
A class $\S$ of $R$-modules is called a Serre subcategory of $\RMod$ if it is closed under submodules, quotients and extensions. 
We also say that a Serre subcategory $\S$ of $\RMod$ is a Serre subcategory of $\Rmod$ if $\S$ consists of finitely generated $R$-modules.

\vspace{5pt}
In \cite{AM-2008}, 
M.\ Aghapournahr and L.\ Melkersson gave the following condition on Serre subcategories of $\RMod$. 

\begin{df}
Let $\S$ be a Serre subcategory of $\RMod$ and $I$ be an ideal of $R$. 
We say that $\S$ satisfies the condition $C_{I}$ if the following condition satisfied:

\vspace{3pt}
\begin{center}
($C_{I}$) \hspace{5pt}  If $M=\G_{I}(M)$ and $(0 :_{M} I)$ is in $\S$, then $M$ is in $\S$. 
\end{center}
\end{df}

\vspace{5pt}
The following special Serre subcategory was named Melkersson subcategory by M.\ Aghapournahr, A.\ J.\ Taherizadeh and A.\ Vahidi in \cite{ATV-2009}. 

\begin{df}
Let $\M$ be a Serre subcategory of $\RMod$. 

\begin{enumerate}
\item\, 
$\M$ is called a Melkersson subcategory with respect to an ideal $I$ of $R$ if $\M$ satisfies the condition $C_{I}$. 

\vspace{3pt}
\item\, 
$\M$ is called a Melkersson subcategory if $\M$ satisfies the condition $C_{I}$ for all ideals $I$ of $R$. 
\end{enumerate}
\end{df}

It has already shown that any Serre subcategory which is closed under injective hulls is the Melkersson subcategory with respect to all ideals $I$ of $R$, 
so that it is a Melkersson subcategory. (See \cite[Lemma 2.2]{AM-2008}.)

%
%
\vspace{7pt}
Next, we consider  classes of extension modules of Serre subcategory by another one.  

\begin{df}
Let $\S_1$ and $\S_2$ be Serre subcategories of $\RMod$. 
We denote by $(\S_1, \S_2)$ the class of all $R$-modules $M$ with 
some $R$-modules $S_{1} \in \S_{1}$ and $S_{2} \in \S_{2}$ such that 
a sequence
$0 \to S_{1} \to M \to S_{2} \to 0$
is exact, that is
\[ 
(\S_{1}, \S_{2}) = \left\{ M \in \RMod \text{ \Large $\mid$ }  
\begin{matrix} \text{there are $S_{1} \in \S_{1}$ and $S_{2} \in \S_{2}$ such that} \cr 
\minCDarrowwidth1pc \begin{CD}0 @>>> S_{1} @>>> M @>>> S_{2} @>>> 0\end{CD} \text{~is exact.}\cr \end{matrix} \right\}. 
\]
We shall refer to $(\S_{1}, \S_{2})$ as a class of extension modules of $\S_{1}$ by $\S_{2}$. 
\end{df}

For example, 
a class $(\FG, \AR)$ is the set of all Minimax $R$-modules 
where  $\FG$ denotes the Serre subcategory consists of all finitely generated $R$-modules 
and $\AR$ denotes the Serre subcategory consists of all Artinian $R$-modules. 
We note that a class $(\S_{1}, \S_{2})$ is not necessarily Serre subcategory.  
(For more detail, see \cite{Y}.)


\vspace{7pt}
\section{Main result}
In this section, 
we shall give an example of Melkersson subcategory which is not closed under injective hulls. 
We denote by $\FS$ the class of $R$-modules with finite support. 
A class $\FS$ is Serre subcategory of $\RMod$ which is closed under injective hulls, 
so that $\FS$ is a Melkersson subcategory. (See \cite[Example 2.4]{AM-2008}.) 
Furthermore, a class $(\FG, \FS)$ is a Serre subcategory of $\RMod$ by \cite[Corollary 4.3 or 4.5]{Y}.

\vspace{5pt}
The main result in this paper is as follows. 

\begin{thm}\label{main}
Let $(R, \m)$ be a local ring with a maximal ideal $\m$.   
Then the following assertions hold. 
\begin{enumerate}
\item\, If $R$ has infinite many prime ideals, then $(\FG, \FS)$ is not closed under injective hulls. 

\vspace{3pt}
\item\, If $R$ is a $2$-dimensional local domain, then $(\FG, \FS)$ is a Melkersson subcategory. 
\end{enumerate}

\noindent
In particular, if $R$ is a $2$-dimensional local domain with infinite many prime ideals,  
then $(\FG, \FS)$ is a Melkersson subcategory which is not closed under injective hulls.
\end{thm}

\begin{proof}
(1)\, 
We assume that $R$ has infinite many prime ideals. 
(We note that the dimension of $R$ must be at least two.) 
Since the set $\Min(R)$ of all minimal prime ideals of $R$ is finite set, 
there exists a prime ideal $\p\in \Min(R)$ such that $V(\p)=\{\q \in \Spec(R) \mid \p \subseteq \q \}$ is infinite set.  
We fix this prime ideal $\p$.

We assume that $(\FG, \FS)$ is closed under injective hulls and shall derived a contradiction. 
Since $R/\p$ is in $(\FG, \FS)$, the injective hull $E_{R}(R/\p)$ of $R/\p$ is also in $(\FG, \FS)$ by assumption.  
Therefore, 
there exists a short exact sequence 
\[ 0 \to F \to E_{R}(R/\p) \to M \to 0 \] 
with $F \in \FG$ and $M \in \FS$. 
Since $V(\p)$ is infinite set and $\Supp(M)$ is finite set, 
we can choose a prime ideal $\n \in V(\p) \setminus \bigl( \Supp(M)\cup \{ \p \} \bigr)$. 
Here, we set $T=R_{\n}$ and $\q=\p R_{\n}=\p T$. 
We note that $T$ is local ring with at least dimension one and $\q$ is a minimal prime ideal of $T$.

Now here, we claim that $E_{T/\q}(T/\q)$ is a finitely generated as $T/\q$-module and shall show this.
By applying the exact functor $(-) \otimes_{R} T $ to the above short exact sequence, 
we see that it holds 
\[ F_{\n} \cong E_{R}(R/\p) \otimes _{R} T \cong  E_{T}(T/\q).\] 
(Also see \cite[Lemma 3.2.5]{BH}.)
Furthermore, it holds  
\[ E_{T/\q}(T/\q) \cong \bigl( 0:_{E_T(T/\q)} \q \bigr) \cong (0:_{F_{\n}} \q)\]
by the above isomorphisms. (Also see \cite[10.1.15 Lemma]{BS}.) 
Since $F$ is a finitely generated $R$-module, $F_{\n}$ is so as $T$-module. 
Thus $E_{T/\q}(T/\q)$ is a finitely generated $T$-module. 
Consequently, we see that $E_{T/\q}(T/\q)$ is a finitely generated as $T/\q$-module.

A local domain $T/\q$ is $\dim\, T/\q \geq 1$ and has a finitely generated injective $T/\q$-module $E_{T/\q}(T/\q)$. 
So it follows from the Bass formula that it holds 
\[ 0<\depth_{T/\q}\, T/\q=\injdim _{T/\q}\, E_{T/\q}(T/\q)=0. \] 
This is a contradiction.

\vspace{3pt}
(2)\,  
We note that any minimal element in $\Supp(M)$ is in $\Ass(M)$ for any (not necessarily finitely generated) $R$-module $M$. 
(e.g. see \cite[Theorem 2.4.12]{E-J}.)

We assume that $R$ is a $2$-dimensional local domain and have to show that a Serre subcategory $(\FG, \FS)$ satisfies the condition $C_{I}$ for all ideals $I$ of $R$. 
We fix an ideal $I$ of $R$. 
We suppose that $X$ is an $R$-module such that $X=\G_{I}(X)$ and $(0 :_{X} I)$ is in $(\FG, \FS)$, 
and shall show that $X$ is in $(\FG, \FS)$.  
There exists a short exact sequence 
\[ 0 \to F \to (0:_{X} I) \to M \to 0 \]
with $F \in \FG$ and $M \in \FS$.

In the case of $\dim\, (0:_{X} I) \leq 1$. 
Then it holds $\Supp(X)=\Ass(X)\cup\{ \m \}$.  
Indeed, since it holds $\Ass(X)=\Ass\bigl((0:_{X}I) \bigr)$, it is easy to see that the zero ideal $(0)$ of $R$ does not belong to $\Supp(X)$. 
Therefore, if there exists a prime ideal $\p \in \Supp(X) \setminus \{ \m \}$, 
$\p$ is minimal in $\Supp(X)$. 
Thus $\p$ is in $\Ass(X)$, so we see that the above equality holds. 
On the other hand, it holds 
\begin{align*}
\Ass(X) 
&= \Ass \bigl( (0:_{X} I) \bigr) \\
&\subseteq \Ass(F) \cup \Ass(M)\\ 
&\subseteq \Ass(F) \cup \Supp(M).
\end{align*}
Since $F$ is a finitely generated $R$-module and $M$ is in $\FS$, $\Ass(X)$ is finite set.  
Consequently, $\Supp(X)$ is also finite set, so we see that $X$ is in $\FS \subseteq (\FG, \FS)$.

In the case of $\dim\, (0:_{X} I) = 2$. 
Since $R$ is a 2-dimensional domain, 
the zero ideal $(0)$ of $R$ must be in $\Supp\bigl( (0:_{X} I) \bigr)$ 
and this is a minimal in $\Supp \bigl( (0:_{X} I) \bigr)$. 
It follows that 
\begin{align*}
(0) \in  \Ass \bigl( (0:_{X} I) \bigr) = V(I) \cap \Ass(X) \subseteq V(I).  
\end{align*}
Therefore, it holds $I=(0)$. 
Consequently, $X=(0:_{X} I)$ is in $(\FG, \FS)$.

The proof is completed.
\end{proof}

\begin{rem}
If $(R, \m)$ is a local ring with at most one dimension, then $\Spec(R)$ is finite set. 
Thus, any support of $R$-module is finite set, so we see $(\FG, \FS)=\RMod$. 
Therefore, in this case, $(\FG, \FS)$ is a Melkersson subcategory and is closed under injective hulls.
\end{rem}

\begin{exm}
Let $R$ be the ring of formal power series $k[[x, y]]$ in the indeterminate $x$ and $y$ with the coefficients in a field $k$. 
Then $R$ is a $2$-dimensional local domain and has infinite many prime ideals $(x+y^n)$ for each non-negative integer $n$. 
Thus, in this case, $(\FG, \FS)$ is a Melkersson subcategory which is not closed under injective hulls by Theorem \ref{main}. 
\end{exm}

\vspace{7pt}
%
%
%
%
\section{Several remarks on Melkersson subcategories}
In this section, we assume that any full subcategory contains a non-zero $R$-module. 

In a local ring $R$, it is clear that any Serre subcategory of $\RMod$ contains all finite length modules. 
On the other hand, 
we can see the following assertion holds.

\begin{prop}
Let $(R, \m)$ be a local ring and $\M$ be a Melkersson subcategory with respect to $\m$. 
Then any Artinian module is in $\M$. 
In particular, Melkersson subcategory contains all Artinian modules. 
\end{prop}

\begin{proof}
Let $\M$ be a Melkersson subcategory with respect to $\m$. 
Since all finite length $R$-modules belong to any Serre subcategory, 
we can see that the injective hull $E_{R}(R/\m)$ of $R/\m$ belongs to $\M$. 
Indeed, since it holds 
\[ \begin{cases}
\ E_{R}(R/\m)=\G_{\m}(E_{R}(R/\m)) \ \ \text{and} \vspace{3pt} \\
\ (0:_{E_{R}(R/\m)} \m) \cong \Hom_{R}(R/\m, E_{R}(R/\m))=R/\m \ \ \text{is in} \ \ \M, 
\end{cases} \]  
it follows from the condition $C_{\m}$ that $E_{R}(R/\m)$ is in $\M$.

Let $M$ be an Artinian module.
Then $M$ is embedded in $\oplus^n E_{R}(R/\m)$ for some integer $n$. 
Therefore, since Melkersson subcategory is closed under finite direct sums and submodules, 
we see that $M$ is in $\M$. 
\end{proof}

\vspace{5pt}
To see  whether Serre subcategory is Melkersson subcategory, 
we have only to check that it satisfies the condition $C_{I}$ for all radical ideals $I$ of $R$. 

\begin{prop}
Let $\M$ be a Serre subcategory. 
Then following conditions are equivalent:

\begin{enumerate}
\item\, $\M$ is a Melkersson subcategory; 

\vspace{3pt}
\item\, $\M$ is a Melkersson subcategory with respect to $\sqrt{I}$ for all ideals $I$ of $R$. 
\end{enumerate}
\end{prop}

\begin{proof}
We assume that $\M$ is a Melkersson subcategory with respect to $\sqrt{I}$ for all ideals $I$ of $R$.
Let $I$ be an ideal of $R$ 
and shall show that $\M$ satisfies condition $C_{I}$.  
We suppose that $M$ is an $R$-module such that $M=\G_{I}(M)$ and $(0:_{M} I)$ is in $\M$. 
Then it holds $\G_{\sqrt{I}}(M) =\G_{I}(M)=M$. 
Furthermore, since $\M$ is closed under submodules and $(0:_{M} \sqrt{I}) \subseteq(0:_{M}I)$, 
we see $(0:_{M} \sqrt{I})$ is in $\M$. 
If follows from the condition $C_{\sqrt{I}}$ that $M$ is in $\M$. 
\end{proof}

\vspace{5pt}
Serre subcategory is defined not only in the category $\RMod$ but also in the category $\Rmod$.
Therefore, 
it stands to reason that we consider the Melkersson subcategory of $\Rmod$ which is defined by considering the condition $C_{I}$ for only finitely generated $R$-modules as follows: the Serre subcategory $\M$ of $\Rmod$ is Melkersson subcategory of $\Rmod$ if it satisfies the condition 

\vspace{3pt}
\begin{center}
($C_{I}$) \hspace{5pt}  If $M=\G_{I}(M) \in \Rmod$ and $(0 :_{M} I)$ is in $\M$, then $M$ is in $\M$ 
\end{center}

\vspace{3pt}
\noindent
for all ideal $I$ of $R$.
However, by the following proposition, 
we can see that it is not necessary to treat Serre subcategory which satisfies such a condition specially.

\begin{prop}
Any Serre subcategory $\S$ of $\Rmod$ is a Melkersson subcategory of $\Rmod$ in the above sense.
\end{prop}

\begin{proof}
By \cite[Theorem 4.1]{Takahashi-2008}, 
there exists a specialization closed subset $W$ of $\Spec(R)$ corresponding to the Serre subcategory $\S$.  
In particular, we can denote
\[ \S =\{  M \in \Rmod \mid \Supp(M) \subseteq W \} \hspace{5pt} \text{and} \hspace{5pt} W =\bigcup_{M \in \S} \Supp (M). \]

\vspace{3pt}
Let $I$ be an ideal of $R$. 
We suppose that $M$ is a finitely generated $R$-module such that $M=\G_{I}(M)$ and $(0 :_{M} I)$ is in $\S$. 
Since $(0 :_{M} I)$ is in $\S$, it holds $\Ass(M)= \Ass((0:_{M} I)) \subseteq \Supp((0:_{M} I)) \subseteq W$, and so we have $\Supp(M) \subseteq W$.  
Consequently, $M$ is in $\S$. 
\end{proof}

\vspace{7pt}

\bigskip \bigskip

\noindent
{\sc Takeshi Yoshizawa}\\
{\sc Graduate School of Natural Science and Technology, \\
Okayama University, Okayama 700-8530, Japan}\\
{\it E-mail address} : \texttt{tyoshiza@math.okayama-u.ac.jp}

\begin{thebibliography}{99}
\addcontentsline{toc}{chapter}{Reference}

%
\bibitem{AM-2008}
{\sc M.\ Aghapournahr} and {\sc L.\ Melkersson},  
Local cohomology and Serre subcategories, 
J.\ Algebra {\bf 320}, 2008, 1275--1287.

%
\bibitem{ATV-2009}
{\sc M.\ Aghapournahr}, {\sc A.\ J.\ Taherizadeh} and {\sc A.\ Vahidi},  
Extension functors of local cohomology modules, 
\texttt{http://arxiv.org/abs/0903.2093v1}.


%
\bibitem{BS}
{\sc M.\ P.\ Brodmann} and {\sc R.\ Y.\ Sharp},
{\it Local cohomology: an algebraic introduction with geometric applications},
Cambridge University Press, Cambridge, 1998.

%
\bibitem{BH}
{\sc W.\ Bruns} and  {\sc J.\ Herzog},
{\it Cohen-Macaulay rings, revised version},
Cambridge University Press, 1998.


%
\bibitem{E-J}
{\sc E.\ E.\ Enochs} and {\sc O.\ M.\ G.\ Jenda}, 
{\it Relative Homological Algebra}, 
Walter De Gruyter, Berlin, New York, 2000. 


%
\bibitem{Takahashi-2008}
{\sc R.\ Takahashi}, 
Classifying subcategories of modules over a commutative noetherian ring, 
J.\ London Math.\  Soc.\ (2) {\bf 78}, 2008, 767--782.

%
\bibitem{Y}
{\sc T.\ Yoshizawa}, 
Classes of extension modules by Serre subcategories, 
\texttt{http://arxiv.org/abs/1011.0376}.


\end{thebibliography}
\end{document}